\documentclass{amsart}
\date{August 16, 2011}
\usepackage[dvips]{graphicx}
\title[CMC-1 Trinoids]{
   CMC-1 trinoids in hyperbolic $3$-space  \\
   and metrics of constant curvature one \\
   with  conical singularities on the $2$-sphere}

\author[Fujimori]{S. Fujimori}
\author[Kawakami]{Y. Kawakami}
\author[Kokubu]{M. Kokubu}
\author[Rossman]{W. Rossman}
\author[Umehara]{M. Umehara}
\author[Yamada]{K. Yamada}

\address[Shoichi Fujimori]{%
   Department of Mathematics,
  Faculty of Science, Okayama University,
  Okayama 700-8530}
\email{fujimori@math.okayama-u.ac.jp}

\address[Yu Kawakami]{
   Graduate School of Science and Engineering,
   Yamaguchi University, Yamaguchi, 753-8512}
\email{y-kwkami@yamaguchi-u.ac.jp}

\address[Masatoshi Kokubu]{
  Department of Mathematics, School of Engineering, 
  Tokyo Denki University,
  Tokyo 101-8457}
\email{kokubu@cck.dendai.ac.jp}

\address[Wayne Rossman]{
   Department of Mathematics,
   Faculty of Science,
   Kobe University,
   Rokko, Kobe 657-8501
}
\email{wayne@math.kobe-u.ac.jp}

\address[Masaaki Umehara]{
   Department of Mathematical and Computing Sciences,
   Tokyo Institute of Technology,
   Tokyo 152-8552}
\email{umehara@is.titech.ac.jp}

\address[Kotaro Yamada]{
   Department of Mathematics,
   Tokyo Institute of Technology,
   Tokyo 152-8551
}
\email{kotaro@math.titech.ac.jp}

\newenvironment{Proof}{\begin{proof}}{\end{proof}}
\newcommand{\op}[1]{{\operatorname{#1}}}
\newcommand{\GL}{\op{GL}}
\newcommand{\PGL}{\op{PGL}}
\newcommand{\SL}{\op{SL}}
\newcommand{\PSL}{\op{PSL}}
\newcommand{\SU}{\op{SU}}

\newcommand{\id}{\op{id}}
\renewcommand{\Im}{\op{Im}}

\newcommand{\pmt}[1]{{\begin{pmatrix} #1  \end{pmatrix}}}
\newcommand{\R}{\boldsymbol{R}}
\newcommand{\C}{\boldsymbol{C}}
\newcommand{\Z}{\boldsymbol{Z}}
\renewcommand{\phi}{\varphi}
\renewcommand{\epsilon}{\varepsilon}

\newcommand\M{{\mathcal M}}
\newcommand{\bP}{\mathbb{P}}
\newcommand{\cmcone}{CMC-$1$}
\newtheorem*{thm}{Theorem}
\newtheorem{cor}{Corollary}
\newtheorem{fact}{Fact}
\theoremstyle{definition}
\newtheorem{defi}{Definition}
\newtheorem{rem}{Remark}
\newtheorem*{ack}{Acknowledgements}
\begin{document}
\maketitle

\begin{abstract}
 \cmcone{} trinoids 
 (i.e.\ constant mean curvature one immersed surfaces
 of genus zero
 with three regular embedded ends)
 in hyperbolic $3$-space $H^3$ are irreducible
 generically, 
 and the irreducible ones 
 have been classified.
 However, the reducible case has not yet been fully 
 treated,
 so here we give an
 explicit description of 
 \cmcone{} trinoids  in $H^3$ that includes 
 the reducible case.
\end{abstract}

\section{Introduction}
Let $H^3$ denote the hyperbolic $3$-space of 
constant sectional curvature $-1$.

A {\em \cmcone{} trinoid\/} in $H^3$ is 
a complete immersed constant mean curvature one surface 
of genus zero with three regular embedded 
ends. 
There are \cmcone{} trinoids with horospherical ends
(i.e.\ regular embedded ends which are asymptotic to
a horosphere).
However, an irreducible trinoid admits only catenoidal 
ends. 
The last two authors \cite{UY6} gave a 
classification of those \cmcone{} trinoids in $H^3$
that are irreducible. 
In particular, they showed that the moduli space of 
irreducible \cmcone{} trinoids in $H^3$
(i.e.\ the quotient space of such immersions by 
the rigid motions of $H^3$) 
corresponds to a certain open dense 
subset of the set of irreducible 
spherical (i.e.\ constant curvature $1$) 
metrics with three conical singularities 
(see Section~\ref{sec:prelim}).  
The paper \cite{UY6} also investigated
the reducible case, 
but had not obtained a complete classification there.

After that, Bobenko, Pavlyukevich, and Springborn 
\cite{BPS} developed a representation formula for 
\cmcone{} surfaces in $H^3$ in terms of holomorphic 
spinors and derived explicit parametrizations 
for irreducible 
\cmcone{} trinoids in $H^3$ in terms of
hypergeometric functions.  
The crucial step in \cite{BPS} was a direct reduction 
of the ordinary differential equation that produces 
\cmcone{} trinoids into a Fuchsian differential equation
with three regular singularities, and we call this 
{\em BPS-reduction}.  
On the other hand, Daniel \cite{D} gave an alternative proof
of the classification theorem for irreducible \cmcone{} 
trinoids, by applying Riemann's classical work on 
minimal surfaces in $\R^3$ bounded by three straight 
lines. 

After the work \cite{UY6} on the irreducible case,
Furuta and Hattori \cite{FH} gave a full classification of
spherical metrics with three conical singularities, 
using a purely geometric method.
Later, Eremenko \cite{E} proved it using
hypergeometric equations.
In this paper, using the argument in \cite{E}
and the BPS-reduction,  
we describe a complete 
classification of reducible \cmcone{} trinoids in $H^3$. 

\section{Preliminaries}
\label{sec:prelim}
Let $M^2$ be a $2$-manifold, and consider a \cmcone{} immersion 
$f:M^2\to H^3$.
The existence of such an immersion implies 
orientability of $M^2$.
By the existence of isothermal coordinates,
there is a unique complex structure on $M^2$
such that the metric $ds^2_f$ induced by $f$
is conformal (i.e.\ $ds^2_f$ is Hermitian).
In this situation, there exists a holomorphic immersion
(called a {\em null lift\/} of $f$)
\[
    F:\tilde M^2 \to \SL(2,\C)
\]
defined on the universal cover $\tilde M^2$ of $M^2$ so that:
\begin{itemize}
 \item $F$ is a {\it null\/} holomorphic map, namely,
       $F_z:=dF/dz$
       is of rank less than $2$  
       on each local complex coordinate $(U;z)$
       of $M^2$.
 \item $f\circ \pi=\hat \pi\circ F$,
       where $\pi:\tilde M^2\to M^2$
       is the covering projection
       and 
       \[
         \hat \pi:\SL(2,\C)\to H^3=\SL(2,\C)/\SU(2)
       \]
       is the canonical projection.
\end{itemize}
Then there exist a meromorphic function $g$ and a holomorphic $1$-form
$\omega$ on $\tilde M^2$ such that
\begin{equation}\label{eq:dF}
 F^{-1}dF=\pmt{g & -g^2 \\ 1 & -g\hphantom{^2}}\omega,
\end{equation}
and the first fundamental form $ds^2_f$ of $f$ satisfies
\[
    ds^2_f=(1+|g|^2)^2|\omega|^2.
\]
The second fundamental form
of $f$ is given by
\[
    h:=-Q-\bar Q+ds^2_f\quad (Q:=\omega dg),
\]
where the holomorphic $2$-differential $Q$ on $M^2$
is called the {\em Hopf differential\/} of $f$.
The set of zeros of $Q$ corresponds to the set 
of umbilics of $f$. 
We set
\begin{equation}\label{eq:Fco}
    F=\pmt{F_{11} & F_{12} \\ F_{21} & F_{22}}.
\end{equation}
Since $\det(dF)=0$, 
one can easily 
show via \eqref{eq:dF} that
\begin{equation}\label{eq:second-gauss}
    g=-\frac{dF_{12}}{dF_{11}}=-\frac{dF_{22}}{dF_{21}}.
\end{equation}
With $\pi_1(M^2)$ denoting the covering transformation
group on the universal cover $\tilde M^2$,
for each $\tau\in \pi_1(M^2)$, there exists
$\rho(\tau)\in \SU(2)$ such that
\begin{equation}\label{eq:repF}
    F\circ \tau=F \rho(\tau),
\end{equation}
which gives a representation
(i.e.\ a group homomorphism)
$\rho:\pi_1(M^2)\to \SU(2)$
satisfying
\begin{equation}\label{eq:tau}
    g\circ \tau^{-1}=\frac{a_{11}g+a_{12}}{a_{21}g+a_{22}}
                 =:\rho(\tau)\star g,
\end{equation}
for each $\tau\in \pi_1(M^2)$,
where $\rho(\tau)=(a_{ij})_{i,j=1,2}$.
\begin{defi} 
 A representation
 $\rho:\pi_1(M^2)\to \SU(2)$
 is called {\em reducible\/}
 if $\rho(\pi_1(M^2))$ is abelian
 and otherwise is called {\em irreducible}.
 A \cmcone{} immersion $f:M^2\to H^3$ is 
 called {\em irreducible\/} ({\em reducible\/})
 if the induced representation
 $\rho$ is irreducible
 (reducible). 
\end{defi}

The meromorphic function (cf.\ \cite{UY})
\[
   G:=\frac{dF_{11}}{dF_{21}}=\frac{dF_{12}}{dF_{22}}
\]
is well-defined on $M^2$, and is  called the 
{\em hyperbolic Gauss map\/} of $f$.

We now consider a \cmcone{} immersion $f$
satisfying the following properties:
\begingroup
\renewcommand{\theenumi}{(\alph{enumi})}
\renewcommand{\labelenumi}{(\alph{enumi})}
\begin{enumerate}
 \item\label{cond:complete-ftc} 
       The metric $ds^2_f$ induced by $f$ is complete 
       and of finite total curvature.
\end{enumerate}
By \ref{cond:complete-ftc}, there exists a 
closed Riemann surface $\bar M^2$ such that
$M^2$ is bi-holomorphic to 
$\bar M^2\setminus \{p_1,\dots,p_n\}$,
where $p_1,\dots,p_n$ are distinct points of
$\bar M^2$ called the {\em ends\/} of $f$.
Then, the Hopf differential $Q$ has at most 
a pole at each of $p_1,\dots,p_n$.
\endgroup

Now, we suppose the second condition:
\begingroup
\renewcommand{\theenumi}{(\alph{enumi})}
\renewcommand{\labelenumi}{(\alph{enumi})}
\begin{enumerate}
\setcounter{enumi}{1}
 \item\label{cond:emb} 
       All the ends $p_1,\dots,p_n$ of $f$ are properly embedded,
      namely, there is a neighborhood $U_j$
      of $p_j$ in $\bar M^2$ such that
      the restriction $f|_{U_j\setminus \{p_j\}}$
      is a proper embedding, for each $j=1,\dots,n$.
\end{enumerate}
\endgroup
\noindent
Then, the condition \ref{cond:emb} implies that
$G$ has at most a pole at each end $p_j$ ($j=1,\dots,n$),
namely, the ends $p_1,\dots,p_n$ are all regular ends.  

\begin{defi}[\cite{T}]
 Let $\bar M^2$ be a closed Riemann surface.
 Let $d\sigma^2$ be a conformal metric on
 $\bar M^2\setminus \{p_1,\dots,p_n\}$,
 where $p_1,\dots,p_n$ are distinct points.
 Then $d\sigma^2$ has a {\em conical singularity\/} 
 of order $\mu_j$ at $p_j$ if
 $\mu_j>-1$ and
 ${d\sigma^2}/{|z|^{2\mu_j}}$
 is positive definite at $p_j$,
 where $z$ is a local coordinate so that $z=0$
 at $p_j$.
 $2\pi(1+\mu_j)$ is called the {\em conical angle\/} 
 of $d\sigma^2$ at $p_j$.
\end{defi}
We set $\bar M^2=S^2$ and consider conformal metrics 
that have exactly three conical singularities at
$0,1,\infty$ on $S^2=\C\cup \{\infty\}$.
We denote by 
$\M_3(S^2)$ the set of such metrics  having
constant curvature $1$ on 
$M^2:=\C\setminus\{0,1\}$,
namely, $\M_3(S^2)$ can be identified with
the moduli space of conformal metrics of constant curvature
$1$ with three conical singularities.
We fix a metric $d\sigma^2\in \M_3(S^2)$,
and then there exists a developing map
\[
   g:\tilde M^2\to S^2=\C\cup \{\infty\}
\]
so that $d\sigma^2=4 dg d\bar g/(1+|g|^2)^2$,
where $\tilde M^2$ is the universal cover 
of $M^2(=\C\setminus\{0,1\})$.
Then there is a representation
(\cite[(2.15) and Lemma 2.2]{UY6})
\begin{equation}\label{eq:rep0}
     \rho:\pi_1(M^2)\to \SU(2)
\end{equation}
satisfying \eqref{eq:tau}.
The metric $d\sigma^2$ is called {\em irreducible\/} if
$\rho$ is irreducible.

We return to the previous situation of \cmcone{} surfaces.
Let $K$ be the Gaussian curvature of the \cmcone{} immersion $f$.
Then
\begin{equation}\label{metric:2}
   d\sigma^2_f:=(-K)ds^2_f
              =\frac{4\,dg\, d\bar g}{(1+|g|^2)^2}.
\end{equation}
This relation implies that $d\sigma^2_f$
has constant curvature $1$ wherever $d\sigma^2_f$ 
is positive definite.
Moreover (\cite{UY}), 
\[
   ds^2_f d\sigma^2_f=4 Q \bar Q
\]
implies that $d\sigma^2_f$ 
has a conical singularity at a zero $q$ of
$Q$, 
and the conical order of $d\sigma^2_f$ 
at $q$ equals $Q$'s order there.
The condition \ref{cond:complete-ftc} implies that
$d\sigma^2_f$ has also a conical 
singularity at each end $p_j$.

\begin{defi}
 Let $f:M^2\to H^3$ be a \cmcone{} immersion satisfying
 conditions \ref{cond:complete-ftc} and \ref{cond:emb}.
 Then $f$ is called a {\em \cmcone{} $n$-noid\/} 
 if $\bar M^2$ is conformally equivalent 
 to the $2$-sphere $S^2$.
 An end $p$ of a \cmcone{} $n$-noid is called {\it catenoidal\/}
 if $Q$ has a pole of order $2$ at $p$.
 A \cmcone{} $n$-noid is called {\it catenoidal\/} if 
 all ends are catenoidal.
\end{defi}
Let $f$ be a \cmcone{} $n$-noid.
When $n=1$, $f$ is congruent to the horosphere.
When $n=2$, $f$ is congruent to a catenoid cousin
or a warped catenoid cousin (cf.\ \cite{RUY2}). 

So it is natural to consider the case $n=3$.
Since the three ends are embedded, the Osserman-type inequality
(\cite{UY}) implies $\deg(G)=2$.
We call a \cmcone{} $3$-noid  a {\em trinoid\/}
(or a \cmcone{} trinoid).
We denote by $\M_3(H^3)$ the set of
congruence classes of trinoids.
We now fix a trinoid $f$.
As shown in \cite{RUY},
there are only two possibilities:
\begingroup
\renewcommand{\theenumi}{(\roman{enumi})}
\renewcommand{\labelenumi}{(\roman{enumi})}
\begin{enumerate}
 \item\label{cond:catenoidal}
       $Q$  has poles of order $2$ at $p_1$, $p_2$, $p_3$.
 \item\label{cond:horospherical}
      $Q$  has at most poles of order $2$ at $p_1$, $p_2$, $p_3$,
      but at least one of the $p_j$ has a pole of order $1$.
\end{enumerate}
\endgroup
As \cmcone{} trinoids satisfying \ref{cond:catenoidal} 
are catenoidal,
irreducible trinoids are catenoidal (see \cite{UY6}).
\cmcone{} immersions satisfying \ref{cond:horospherical}
have been classified in \cite[Theorems 4.5--4.7]{RUY}.
So from now on we consider just the case \ref{cond:catenoidal}.
Without loss of generality we may assume
$p_1=0$, $p_2=1$, $p_3=\infty$. 
As mentioned above, the metric
$d\sigma^2_f$ given 
by \eqref{metric:2} has conical singularities at
the zeros of $Q$ and the three ends $p_1$, $p_2$, $p_3$.  
We denote by $\beta_j(>-1)$ the order of  $d\sigma^2_f$ at $p_j$,
and by 
\[
    B_j:=\pi(1+\beta_j)(>0)\qquad (j=1,2,3)
\]
the half of the conical angle of  $d\sigma^2_f$ at $p_j$ ($j=1,2,3$).
The group $\rho(\pi_1(M^2))$ is generated by three monodromy matrices 
$\rho_1$, $\rho_2$, $\rho_3$ which represent
loops surrounding $z=0$, $1$, $\infty$.
Each $\rho_j$ ($j=1,2,3$) has eigenvalues $-\exp(\pm iB_j)$.
Then we have  (cf.\ \cite{UY6})
\begin{equation}\label{eq:Q}
   2Q=
     \frac{c_3z^2+(c_2-c_3-c_1)z+c_1}{z^2(z-1)^2}dz^2,
\end{equation}
where $c_j:=-{\beta_j(\beta_j+2)}/{2}$ 
does not vanish by \ref{cond:catenoidal} (i.e.\ $B_j\ne \pi$)
for $j=1$, $2$, $3$, and
\begin{equation}\label{eq:hanbetu}
   \frac{(c_1)^2+(c_2)^2+(c_3)^2}2\ne c_1c_2+c_2c_3+c_3c_1.
\end{equation}
We denote by $q_1$, $q_2$ the two 
roots of the equation
\begin{equation}\label{eq:2}
    c_3z^2+(c_2-c_3-c_1)z+c_1=0.
\end{equation}
Since $c_3\ne 0$,
the Hopf differential $Q$ has
exactly two zeros at $q_1$ and $q_2$.
In fact,
\eqref{eq:hanbetu} is equivalent to the condition
$q_1\ne q_2$ 
(i.e.\ the discriminant of \eqref{eq:2} does not vanish).
As shown in \cite{UY6}, the condition \ref{cond:emb}
implies that $G$ does not branch at 
the three ends $0$, $1$, $\infty$, but has exactly
two branch points at $q_1$, $q_2$.
Since $G$ is of degree $2$
and has the ambiguity of  M\"obius transformations,
we may set (cf.\ \cite{UY6}) 
\begin{equation}\label{eq:G}
    G:=z+\frac{(q_1-q_2)^2}{2(2z-q_1-q_2)}.
\end{equation}
Take a solution $F:\tilde M^2\to \SL(2,\C)$ of
the ordinary differential equation
\begin{equation}\label{eq:dual}
   dF F^{-1}=\pmt{G & -G^2 \\ 1 & -G\hphantom{^2}}\frac{Q}{dG}.
\end{equation}
If the image $\rho(\pi_1(M^2))$ of
the representation $\rho$ of $F$
is conjugate to a subgroup of
$\SU(2)$, then $f=\hat \pi(Fa)$ gives a \cmcone{} trinoid
for a suitable choice of $a\in\SL(2,\C)$
(cf.\ \eqref{eq:repF}).
We denote by $\M_{B_1, B_2, B_3}(H^3)$
(resp.\ $\M_{B_1, B_2, B_3}(S^2)$)
the congruence classes of trinoids $f$
satisfying \ref{cond:catenoidal}
(resp.\ of the metrics $d\sigma^2$ of constant curvature $1$)
such that $d\sigma^2_f$ (resp. $d\sigma^2$)
has conical angle $2B_j (\ne 2\pi)$ at each $p_j$.

\begin{fact}[\cite{UY6}]
\label{fact:1}
 For each $B_1$, $B_2$, $B_3\in (0,\infty)$,
 $\M_{B_1, B_2, B_3}(H^3)$
 {\rm(}resp.\ $\M_{B_1, B_2, B_3}(S^2)${\rm)}
 consists of a unique irreducible element
 if it satisfies \eqref{eq:hanbetu}
 {\rm(}resp.\ no condition{\rm)} and
 \begin{multline}\label{eq:irr}
  \cos^2 B_1+\cos^2 B_2+\cos^2 B_3\\
               +2\cos B_1\cos B_2\cos B_3<1.
 \end{multline}
 Conversely, any irreducible 
 trinoids {\rm(}resp.\ any irreducible metrics in
 $\M_{3}(S^2)${\rm)}
 are so obtained.

 In particular, there is a unique catenoidal trinoid $f$
 such that
 \begin{itemize}
  \item the hyperbolic Gauss map
	$G$ is given by \eqref{eq:G},
  \item the Hopf differential
	$Q$ is given by \eqref{eq:Q},
  \item $d\sigma^2_f$
	has conical angle $2B_j$ at each
	end $p_j$.
 \end{itemize}
\end{fact}

\begin{figure}[hbt]
 \begin{center}
  \begin{tabular}{c@{\hspace{2em}}c}
       \includegraphics[width=3cm]{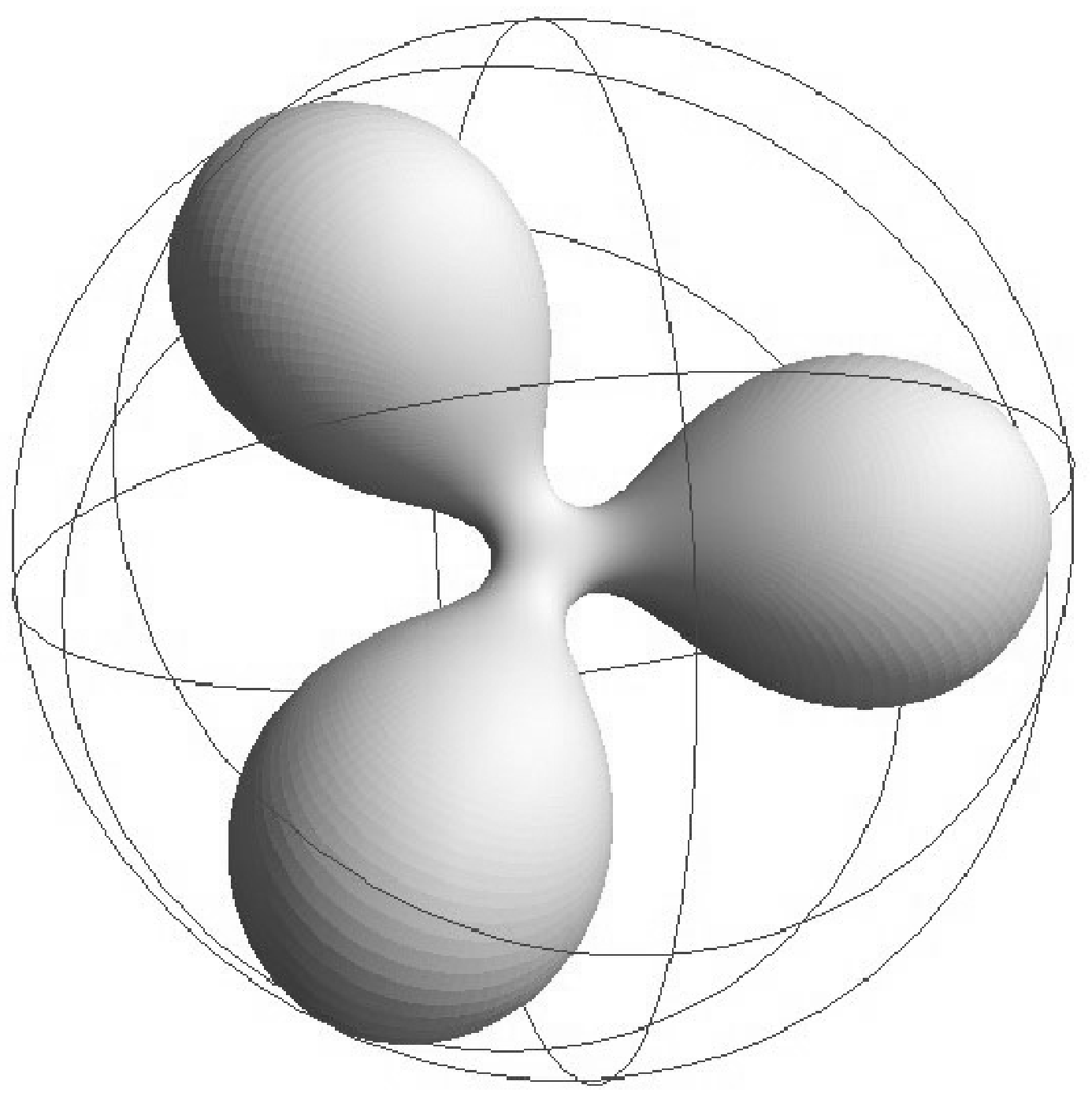} & 
       \includegraphics[width=3cm]{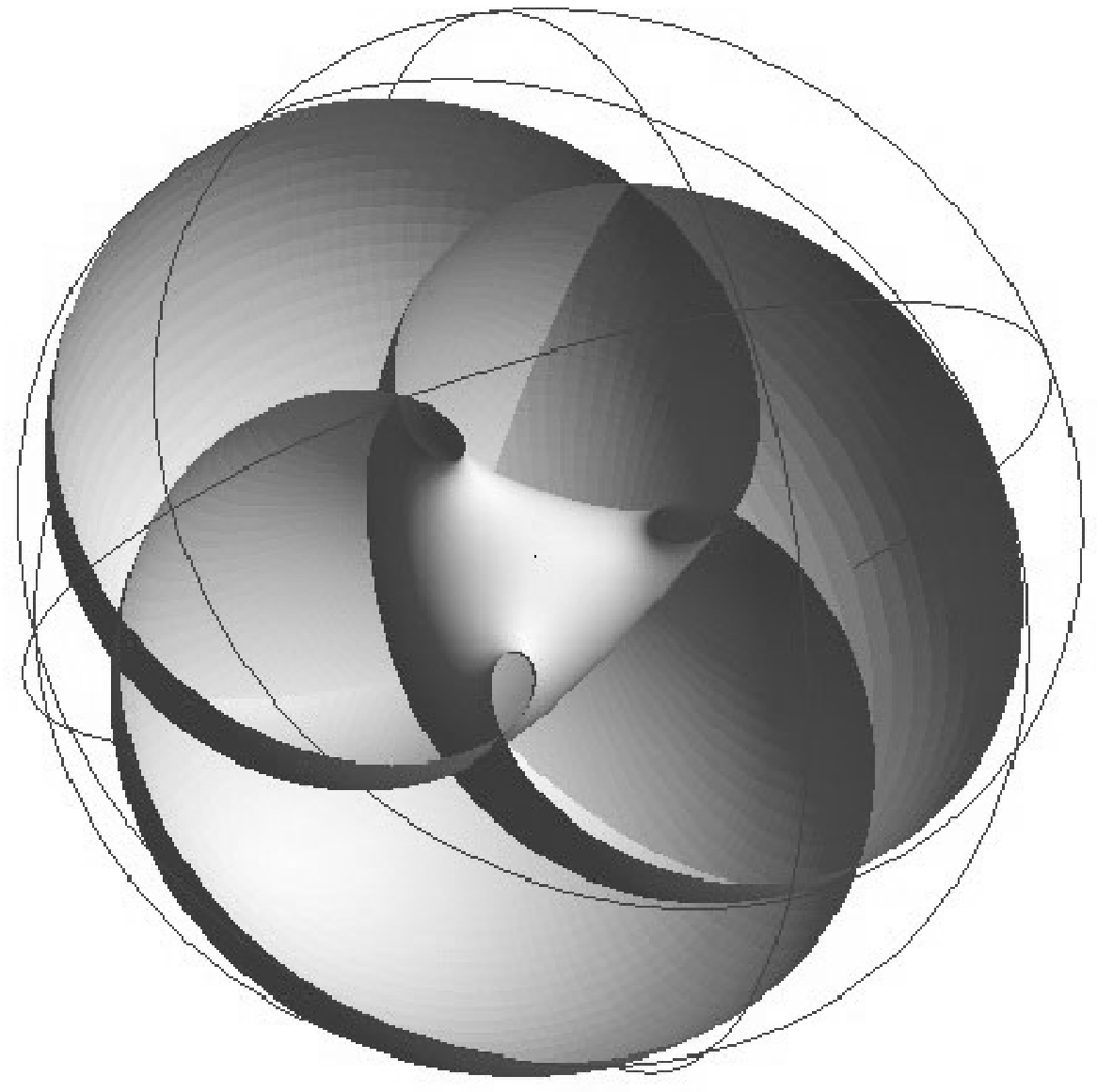}
  \end{tabular}
 \end{center}
\caption{%
 Trinoids with $B_1=B_2=B_3$.
}
\label{fig:trinoid2}
\end{figure}

Figure \ref{fig:trinoid2}, left (resp.\ right)  is  
an irreducible trinoid 
(resp.\ a cutaway view of an irreducible trinoid)
for $B_1=B_2=B_3(=B)$ with $B<\pi$ (resp.\ $B>\pi$).  

\begin{rem}
 Since the hyperbolic Gauss map $G$ changes
 under a rigid motion of $H^3$, 
 the above trinoid $f$ is uniquely determined without the ambiguity of
 isometries of $H^3$ (cf.\ \cite[Appendix B]{UY6}).
 After \cite{UY6},
 Bobenko, Pavlyukevich and Springborn
 \cite{BPS} gave a different proof, 
 whose underlying idea also appears in the next section.
 Also,  Daniel \cite{D} gave an alternative proof of
 this fact (see the introduction). 
\end{rem}

For each $B_j$ ($j=1,2,3$)
there exists a unique real number $\hat B_j\in [0,\pi]$
such that $\cos B_j=\cos \hat B_j$,
since $\cos t=\cos(2\pi-t)$ for $t\in [0,2\pi)$.
By definition, it holds that
$B_j\ge \hat B_j$.
Without loss of generality, 
we may
assume that
$\hat B_1\le \hat B_2\le \hat B_3$.
We now set
$B'_1:=\hat B_1$,
and for $j=2,3$,
\[
  B'_j:=
  \begin{cases}
     \hat B_j &     \quad\mbox{if $\hat B_2+\hat B_3\le \pi$}, \\
     \pi-\hat B_j & \quad \mbox{if $\hat B_2+\hat B_3> \pi$}. 
  \end{cases}
\]
Then we have that
\begin{equation}\label{eq:er2}
    0\le B'_1+B'_2,~
    B'_1+B'_3,~
    B'_2+B'_3\le \pi,
\end{equation}
and the condition
\eqref{eq:irr} is equivalent to
{\small
\[
 \cos^2 B'_1+\cos^2 B'_2+\cos^2 B'_3
     +2\cos B'_1\cos B'_2\cos B'_3<1,
\]}%
which is equivalent to the condition
\begin{multline*}
 \cos\frac{B'_1+B'_2+B'_3}2
 \cos\frac{-B'_1+B'_2+B'_3}2 \\
    \times \cos\frac{B'_1-B'_2+ B'_3}{2}
           \cos\frac{B'_1+B'_2-B'_3}2<0.
\end{multline*}
By \eqref{eq:er2}, this then reduces to the condition
\begin{equation}\label{eq:irr2}
    B'_1+B'_2+B'_3>\pi.
\end{equation}
The condition \eqref{eq:irr}
(or equivalently \eqref{eq:irr2}) 
implies $B_j\not\in \pi \Z$ ($j=1,2,3$), and
is the same condition 
as in \cite{UY6}, \cite{FH} 
or \cite{E} that there exists an irreducible metric in 
$\M_3(S^2)$ with three conical angles $2B_1$, $2B_2$, $2B_3$.

\section{Reducible trinoids}
\label{sec:red}
Let $\alpha$ be a $2\times 2$-matrix valued
meromorphic $1$-form on $\C\cup \{\infty\}$.
Consider an ordinary differential equation
\begin{equation}\label{eq:ode}
     dE E^{-1}=\alpha,
\end{equation}
which is called a {\em Fuchsian differential equation\/}
if it admits only regular singularities.
For example, the equation \eqref{eq:dual} with $G$, $Q$ 
satisfying \eqref{eq:G} and \eqref{eq:Q} is a Fuchsian differential
equation with regular singularities at 
$z=0$, $1$, $\infty$, $q_1$, $q_2$. 
Let
$p_1,\dots,p_n\in \C\cup \{\infty\}$
be the regular singularities of the equation \eqref{eq:ode}.
We denote by $\tilde M^2$ the universal cover of 
\[
   M^2:=\C\cup \{\infty\}\setminus \{p_1,\dots,p_n\}.
\]
Then there exists a solution 
$E:\tilde M^2\to \GL(2,\C)$ of \eqref{eq:ode}.
Since $\alpha$ is defined on $M^2$,
there exists a representation
$\gamma:\pi_1(M^2)\to \GL(2,\C)$
such that
$E\circ \tau=E\gamma(\tau)$.
Let 
\[
    \GL(2,\C)\ni a \mapsto [a]\in \PGL(2,\C)
      =\PSL(2,\C)
\]
be the canonical projection.
Then 
\[
   h_1:=-E_{12}/E_{11},\qquad 
   h_2:=-E_{22}/E_{21}
\]
satisfy (see \eqref{eq:tau} for the definition of $\star$)
\[
    h_i\circ \tau^{-1}=\gamma(\tau)\star h_i
               \qquad (\tau \in \pi_1(M^2),\,\, i=1,2),
\]
where $E=(E_{jk})_{j,k=1,2}$.
Thus the functions $h_i$ ($i=1,2$)
induce a common group homomorphism
$[\gamma]:\pi_1(M^2)\to \PGL(2,\C)$
which is called the {\em monodromy representation\/} 
of the equation \eqref{eq:ode}.
In particular, 
the representation $[\rho]$ for $F$
as in \eqref{eq:dual} is just the
monodromy representation.
\begin{defi}
  Let $r(z),s(z)$ be meromorphic functions
  on $\C\cup \{\infty\}$ and
  \begin{equation}\label{eq:E2}
     X''+rX'+sX=0
  \end{equation}
 be an ordinary differential equation
 with regular singularities at $z=p_1,\dots,p_n$, 
 where $X'=dX/dz$.
 Then there exists a pair of
 solutions $w_1,w_2:\tilde M^2\to \C$
 which are linearly independent, and
 $\{w_1,w_2\}$ is called a fundamental
 system of solutions.
 There exists a representation
 $\gamma:\pi_1(M^2)\to \GL(2,\C)$
 for each fundamental system $\{w_1,w_2\}$, 
 such that
 \[
    (w_1\circ \tau,w_2\circ \tau)=(w_1,w_2)\gamma(\tau),
 \]
 where $(w_1,w_2)$ is a row vector. 
 As a monodromy of the function $-w_2/w_1$,
 the induced homomorphism
 $[\gamma]:\pi_1(M^2)\to \PGL(2,\C)$
 is called the {\em monodromy representation\/} of
 the equation \eqref{eq:E2}.
\end{defi}

To give a complete classification of trinoids,
the following reduction given in \cite{BPS}
is crucial:
Let $F$ be a null lift of the catenoidal trinoid $f$
whose hyperbolic Gauss map $G$ 
and Hopf differential $Q$ are
given by \eqref{eq:G} and \eqref{eq:Q},
respectively.
In the expression \eqref{eq:dual},
we can write
\[
   \pmt{G & -G^2 \\ 1 & -G\hphantom{^2}}\frac{Q}{dG}=
      \pmt{\bP_1 \bP_2 & (\bP_1)^2 \\ 
            -(\bP_2)^2 &  -\bP_1 \bP_2}dz,
\]
where
$\bP_i:=\dfrac{p^0_{i}}{z}+\dfrac{p^1_{i}}{z-1}+
           {p^\infty_{i}}$
and $p^0_{i},p^1_{i},p^\infty_{i}$ ($i=1,2$)
are constants depending only on $B_1$, $B_2$, $B_3$.
In \cite{BPS}, 
the matrix $\Phi:=D^{-1}F$
is defined by
\[
    D:=\sqrt{z-1}
         \pmt{\bP_1 & \alpha_1z+\beta_1 \\
             -\bP_2 & \alpha_2z+\beta_2}
         \pmt{1 & 0 \\ \frac{k}{z(z-1)} & \frac{1}{z-1}}
         \pmt{\vartheta & 0 \\ 1 & 1},
\]
where $\alpha_j$, $\beta_j$ ($j=1,2$), 
$k$ and $\vartheta$ are all real constants
depending only on  $B_1$, $B_2$, $B_3$.
Then there exist $2\times 2$
matrices $A_0$, $A_1$ with real coefficients such that
\begin{equation}\label{eq:red}
   d\Phi \Phi^{-1}
          =\left(\frac{A_0}{z}+\frac{A_1}{z-1}\right)dz.
\end{equation}
We call \eqref{eq:red}
the {\em BPS-reduction\/} of \eqref{eq:dual}.
(This reduction does not work if $f$ has a horospherical end, 
 but such trinoids would be in the case 
 \ref{cond:horospherical} mentioned before.) 
By \eqref{eq:repF}, it holds that
$\Phi\circ \tau=\Phi \rho(\tau)$ for each $\tau \in \pi_1(M^2)$.
Obviously \eqref{eq:red} has three
regular singularities at $z=0$, $1$, $\infty$.
Since $A_0$ and $A_1$ are both constant real
matrices, it is well-known that there exist
real numbers $a$, $b$, $c$ such that
the monodromy representation of
the ordinary differential equation
(called the {\em hypergeometric equation})
\begin{equation}\label{eq:gauss}
         z(1-z)X''+(c-(a+b+1)z)X'-ab X=0
\end{equation}
is conjugate to that of \eqref{eq:red} (i.e.\ $[\rho]$).
On the other hand, if we express $F$ as in \eqref{eq:Fco},
then $X=F_{11}$, $F_{12}$ satisfy the ordinary differential 
equation (cf.\ \cite[p.\ 32]{RUY})
\begin{equation}\label{eq:X}
   X''-\bigl(\log (\hat Q/G')\bigr)'X'+\hat Q X=0,
\end{equation}
where $Q=\hat Q(z)dz^2$ and $G'=dG/dz$.
Thus the monodromy representation of \eqref{eq:X}
with respect to $(F_{11},F_{12})$ is equal to that of $F$.
In particular, the monodromy representation of \eqref{eq:X} 
is conjugate to that of \eqref{eq:gauss}. 
Hence, these two ordinary differential 
equations have the same 
{\em exponent\/} (i.e.\ the difference  of the two 
solutions of the indicial equation)
at each regular singularity.
Since \eqref{eq:X} has the exponent
$B_1/\pi$, $B_2/\pi$, $B_3/\pi$ at $z=0,1,\infty$, 
respectively, we have
\begin{align*}
   \pm B_1&=\pi(1-c), \quad
   \pm B_2= \pi(a-b), \\
   \pm B_3&= \pi(c-a-b),
\end{align*}
which is the same set of relations as in \cite[(4)]{E}.
This implies the classification of catenoidal trinoids reduces 
to that of metrics in $\M_3(S^2)$.
In particular,
the classification 
results for reducible metrics 
in $\M_3(S^2)$
given in Furuta-Hattori \cite{FH}
and Eremenko \cite[Theorem 2]{E} yield the following 
assertion.

\begin{thm}\label{thm}
 Suppose $B_1/\pi$ is an integer, 
 and $B_j\ne \pi$ {\rm(}$j=1,2,3${\rm)}. 
 Then 
 $\M_{B_1,B_2,B_3}(H^3)$ {\rm(}resp.\ $\M_{B_1,B_2,B_3}(S^2)${\rm)}
 is non-empty if and only if $B_1$, $B_2$, $B_3$ satisfy 
\eqref{eq:hanbetu} 
 {\rm(}resp.\ no condition{\rm)} and
 one of the following two conditions{\rm:}
 \begingroup
   \renewcommand{\theenumi}{($\mathrm{C}_{\arabic{enumi}}$)}
   \renewcommand{\labelenumi}{($\mathrm{C}_{\arabic{enumi}}$)}
 \begin{enumerate}
 \item\label{cond:h-one}
      $B_2,B_3\not\in \pi\Z$, but
      either $|B_2-B_3|/\pi$ or $(B_2+B_3)/\pi$
      is an integer $m$ of opposite parity 
      from $B_1/\pi$, and $\pi m\le B_1-\pi$.
      In this case, 
      $\M_{B_1,B_2,B_3}(H^3)$ {\rm(}resp.\ $\M_{B_1,B_2,B_3}(S^2)${\rm)}
      is $1$-dimensional.
 \item\label{cond:h-three}
      $B_2,B_3\in \pi\Z$,
      and $(B_1+B_2+B_3)/\pi$ is odd, 
      and each of $B_1$, $B_2$, $B_3$ is less than the sum of
      the others.
      In this case, 
      $\M_{B_1,B_2,B_3}(H^3)$ {\rm(}resp.\ $\M_{B_1,B_2,B_3}(S^2)${\rm)}
      is $3$-dimensional.
 \end{enumerate}
\endgroup
\end{thm}
\begin{cor}\label{cor2}
 A catenoidal trinoid $f$ is irreducible if and only if 
 $B_1/\pi$, $B_2/\pi$, $B_3/\pi$
 are all non-integers, and
 $f$ is reducible
 if and only if at least one of $B_1/\pi$, $B_2/\pi$, $B_3/\pi$
 is an integer.
\end{cor}
\begin{Proof}
 A trinoid $f$ is irreducible if the representation $\rho$
 as in \eqref{eq:repF} is irreducible.
 The representation $\rho$ coincides with that of
 the corresponding metric in $\M_3(S^2)$.
 The corresponding  assertion for metrics 
 in $\M_3(S^2)$ is proved in \cite[Lemma 3.1]{UY6}.
\end{Proof}

\begin{rem}
 Reducibility is equivalent to 
 at least one of $B_1/\pi$, $B_2/\pi$, $B_3/\pi$ being an integer. 
 This cannot be proved purely algebraically,
 as there are diagonal matrices
 $\rho_1$, $\rho_2$, $\rho_3$ in $\SU(2)$ with 
 $\rho_1\rho_2\rho_3=\id$ so that
 no eigenvalues of $\rho_1$, $\rho_2$, $\rho_3$ are 
 $\pm 1$.
\end{rem}

\begin{figure}[hbt]
\begin{center}
  \includegraphics[width=0.9in]{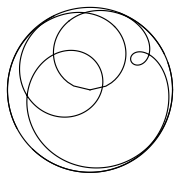}
\end{center}
\caption{
 A profile curve of a non-symmetric trinoid with 
 $B_1=B_2=B_3=3\pi$. (The outer circle represents the 
 ideal boundary of $H^3$.)}
\label{fig:trinoid3}
\end{figure}

\begin{rem}\label{rm:ER}
 Eremenko \cite[Theorem 2]{E} asserts the uniqueness of 
 $d\sigma^2 \in \M_{B_1, B_2, B_3}(S^2)$ 
with prescribed conical angles.  
 This is 
 correct in the irreducible case, but if 
 $B_1\in \pi\Z$, then the 
 metric has a nontrivial deformation
 preserving its conical angles:
 A metric $d\tau^2\in \M_{3}(S^2)$ 
 has the same conical angles as those of $d\sigma^2$ 
 if and only if each developing map 
 of $d\tau^2$ is given by $k=a\star h$ for suitable
 $a\in \SL(2,\C)$, where $h$ is a
 developing map of $d\sigma^2$.
 So $\M_{B_1, B_2, B_3}(S^2)$
 can be identified with the set (\cite{UY6})
 \[
    \{\hat \pi (a)\,;\, 
     a(\Im\rho) a^{-1}\subset \SU(2),~a\in \SL(2,\C)\}
    (\subset H^3), 
 \]
 where $\hat \pi:\SL(2,\C)\to H^3$ is the canonical
 projection and $\Im\rho$ is the image of 
 $\rho$ as in \eqref{eq:rep0}.
 Then $\M_{B_1, B_2, B_3}(S^2)=H^3$ 
 if $B_j/\pi$ ($j=1,2,3$) are all integers,
 and $\M_{B_1, B_2, B_3}(S^2)$ is a geodesic line in $H^3$
 if
 one of $B_j/\pi$ ($j=1,2,3$) is not an integer (cf.\ \cite{UY6}).
 A metric $d\sigma^2$ in $\M_{B_1, B_2, B_3}(S^2)$
 is called {\it symmetric\/} if
 the metric is invariant under 
 an anti-holomorphic involution.
 We denote by $\hat {\M}_{B_1, B_2, B_3}(S^2)$
 the subset consisting of symmetric metrics 
 in $\M_{B_1, B_2, B_3}(S^2)$.
 If $B_j/\pi$ ($j=1,2,3$) are all integers,
 $\hat {\M}_{B_1, B_2, B_3}(S^2)$ consists of a
 hyperbolic plane in $H^3$ (\cite{UY6}).
 If one of $B_j/\pi$ ($j=1,2,3$) is not an integer,
 $\hat {\M}_{B_1, B_2, B_3}(S^2)$ coincides with
 $\M_{B_1, B_2, B_3}(S^2)$.
 A metric $d\sigma^2$ 
 in $\hat {\M}_{B_1, B_2, B_3}(S^2)$
 with conical angles $2B_1$, $2B_2$ and $2B_3$
 can be regarded as a doubling of
 the generalized spherical triangle with angles $B_1$, $B_2$
 and $B_3$.
 Using this, Furuta-Hattori \cite{FH} gave two operations in 
 $\hat {\M}_3(S^2)$ for distinct $\{i,j,k\}=\{1,2,3\}$: 
 \begin{align*}
    (B_i,B_j,B_k) \mapsto (B_i+\pi,B_j+\pi,B_k) , \\
    (B_i,B_j,B_k) \mapsto (\pi-B_i,B_j+\pi,B_k) , 
 \end{align*}
 with the second operation allowed
 only when $B_i<\pi$.
 The first operation is attaching a closed hemisphere
 in $S^2$ to the edge $B_iB_j$ of the 
 spherical triangle $\triangle B_iB_jB_k$.
 The second operation is attaching a geodesic
 $2$-gon of equi-angles $\pi-B_i$ to the 
 edge $B_iB_j$ so that the initial vertex $B_i$ becomes 
 an interior point of an edge of the new triangle.
 Conditions \ref{cond:h-one} and \ref{cond:h-three} are invariant under
 these two operations.
 Moreover, the three angles $(B_1,B_2,B_3)$
 satisfying conditions \ref{cond:h-one} and \ref{cond:h-three}
 are obtained from a given initial data
 $(B'_1,B'_2,B'_3)$ by these two operations.
 Furuta-Hattori proved this using a geometric argument.
 On the other hand, Eremenko found
 \ref{cond:h-one} and \ref{cond:h-three}
 from the viewpoint of hypergeometric equations.
 We remark that spherical triangles of arbitrary
 angles $B_1$, $B_2$, $B_3\in (0,\infty)$
 were investigated  by Felix Klein in 1933
 (see the end  of \cite{UY6}).
 The trinoid shown in Figure \ref{fig:trinoid3} 
 is not symmetric, although there does
 exist a symmetric trinoid with the same conical
 angles and dihedral symmetry.
\end{rem}

\begin{figure}[hbt]
\begin{center}
\footnotesize
\begin{tabular}{cc}
  \includegraphics[width=0.7in]{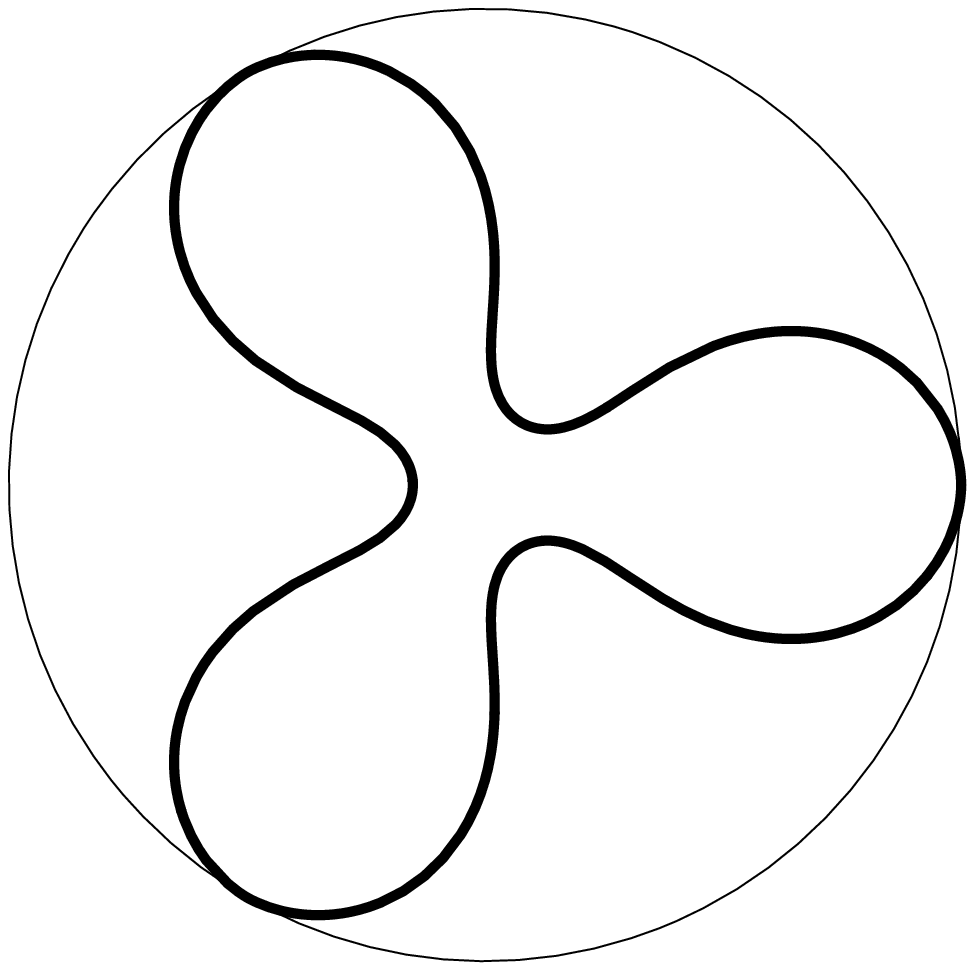} &
  \includegraphics[width=0.7in]{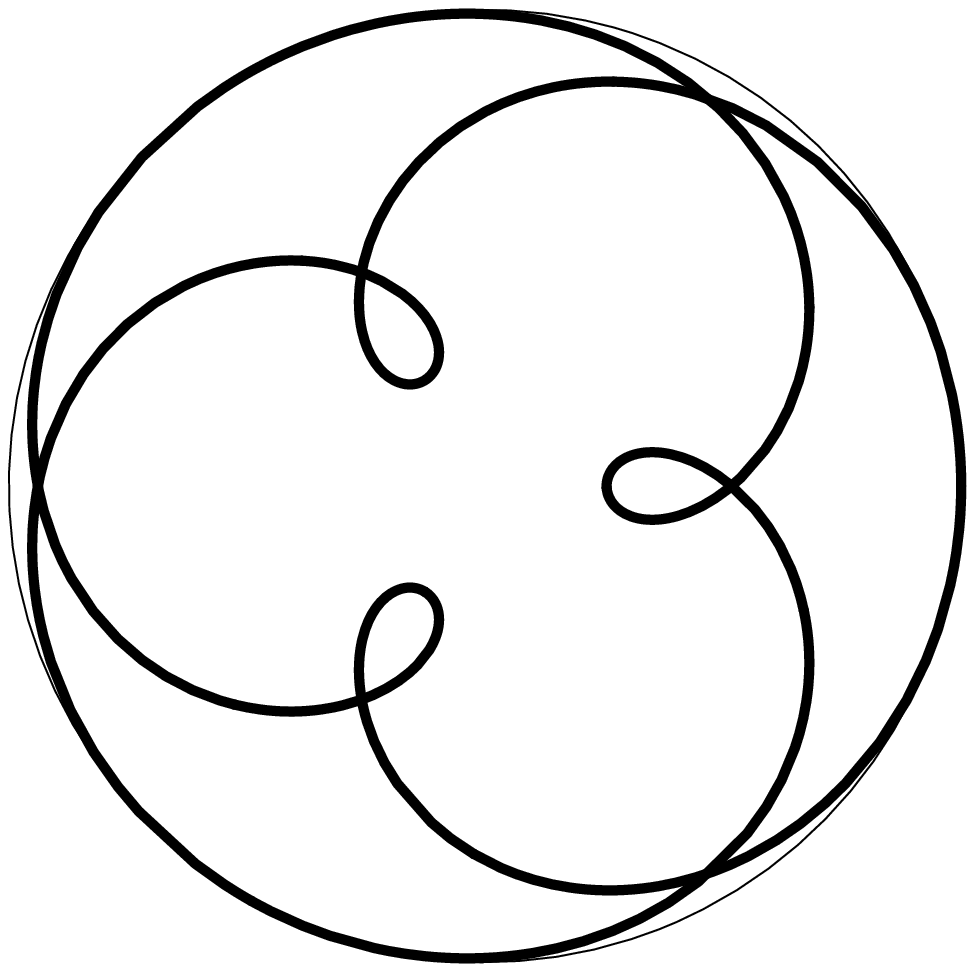} \\
  \includegraphics[width=0.7in]{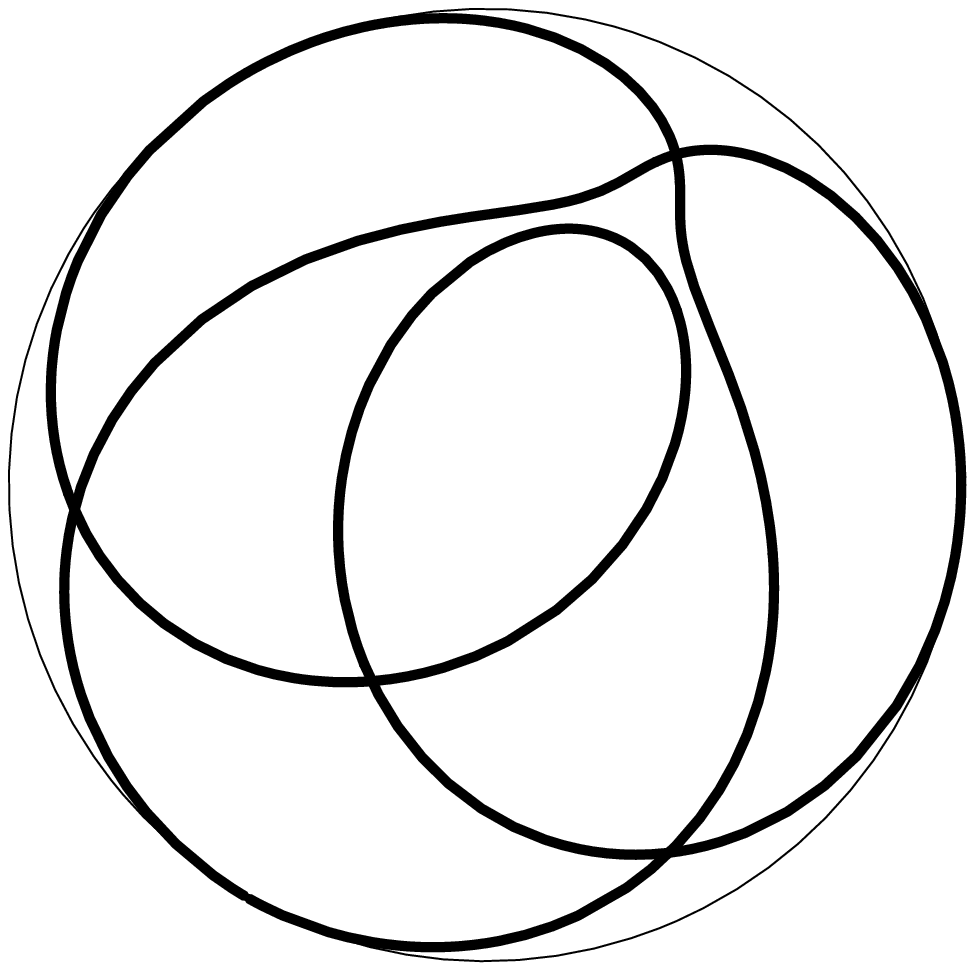} &
  \includegraphics[width=0.7in]{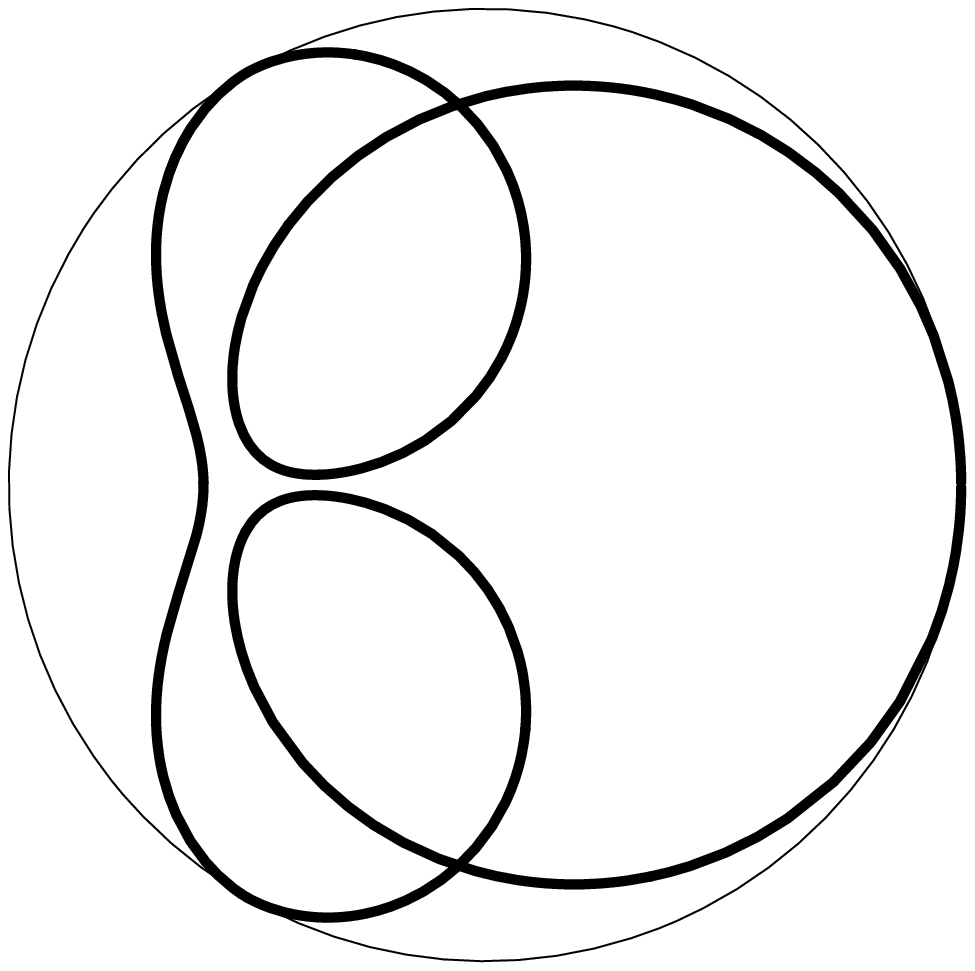} \\
  Type $(+,+,+)$ &
  Type $(-,-,-)$ \\
  Type $(+,-,-)$ &
  Type $(+,+,-)$
\end{tabular}
\end{center}
 \caption{Profile curves of trinoids of different types.}
\label{fig:trinoid}
\end{figure}

Finally, we group the surfaces by the 
signatures of $c_1$, $c_2$, $c_3$.
For example, a trinoid $f$ is 
said to be of type $(+,+,+)$ if $c_1$, $c_2$, $c_3$ are all 
positive, and of type $(-,+,+)$ 
if one of $c_1$, $c_2$, $c_3$
is negative and the other two are positive, etc.  
As remarked in \cite{RUY2},
by numerical experiment, 
it seems that the four types
$(+,+,+)$, $(-,+,+)$,
$(-,-,+)$ and $(-,-,-)$ 
have distinct regular homotopy types
(see Figure~\ref{fig:trinoid}).
 Surfaces of type $(+,+,+)$ have absolute total 
curvature less than $8\pi$, and it 
seems that only surfaces in this class 
can be embedded.
\begin{ack}
 The authors thank Masaaki Yoshida and Yoshishige Haraoka
 for valuable comments, and for imparting to us the charm of
 the hypergeometric equation. 
 They also  thank Alexander Eremenko
 for valuable comments.
\end{ack}

\end{document}